\documentclass[12pt]{amsart}
\usepackage{amsfonts}
\usepackage{amssymb}
\usepackage{graphicx}
\usepackage{fourier}
\theoremstyle{definition}

\theoremstyle{remark}

\newcommand{\R}{\mathbb R}

\newcommand{\ds}{\displaystyle}
\begin{document}
\title[WEIGHTED SUMS OF THE SQUARES OF THE DISTANCES]
{WEIGHTED SUMS OF THE SQUARES OF THE DISTANCES OF A POINT TO THE SIDELINES OF A TRIANGLE}

\author{GEORGI GANCHEV AND NIKOLAI NIKOLOV}%
\address{Bulgarian Academy of Sciences, Institute of Mathematics and Informatics,
Acad. G. Bonchev Str. bl. 8, 1113 Sofia, Bulgaria}
\email{ganchev@math.bas.bg}
\email{nik@math.bas.bg}

\begin{abstract}
We study a function, which is a weighted sum of the squares of the distances of an arbitrary point to the sidelines of a triangle.
The given weights, considered as barycentric coordinates, determine a point $M$. We prove that the function reaches
its minimum (maximum) at a point, which is isogonal conjugate to $M$.
\end{abstract}

\maketitle

\section{Introduction}
As usual, we denote by $a, b, c$ the sides of a given $\triangle ABC$ and by $S$ its area.
The positive orientation of the plane is determined by $\triangle ABC$.

Let $x, y, z$ be trilinear coordinates of a point with respect to $\triangle ABC$.

It is well known that the Lemoine point $K$ minimizes the sum $x^2+y^2+z^2$, i.e. \cite{L}
$$x^2+y^2+z^2 \geq \frac{4S^2}{a^2+b^2+c^2}\, .$$

Recently Kimberling \cite{K} obtained several inequalities for the power sums $x^q+y^q+z^q$.

In this note for an arbitrary point $X$ in the plane of $\triangle ABC$ we study a weighted sum $F(X)$
of the squares of the distances of $X$ to the sidelines of the triangle. We give a
geometric interpretation of the minimum (maximum) of the function $F(X)$.

\section{Preliminaries}
Let us recall some properties of isogonal conjugate points with respect to a given $\triangle ABC$.

Given the basic $\triangle ABC$ and its circumcircle $k(ABC)$. Denote by
$\imath$ the isogonal conjugation with respect to the triangle. The action of $\imath$ in the domains
(with respect to the vertex $A$) (Fig. 1) is as follows:

\begin{center}
\includegraphics[width=10cm]{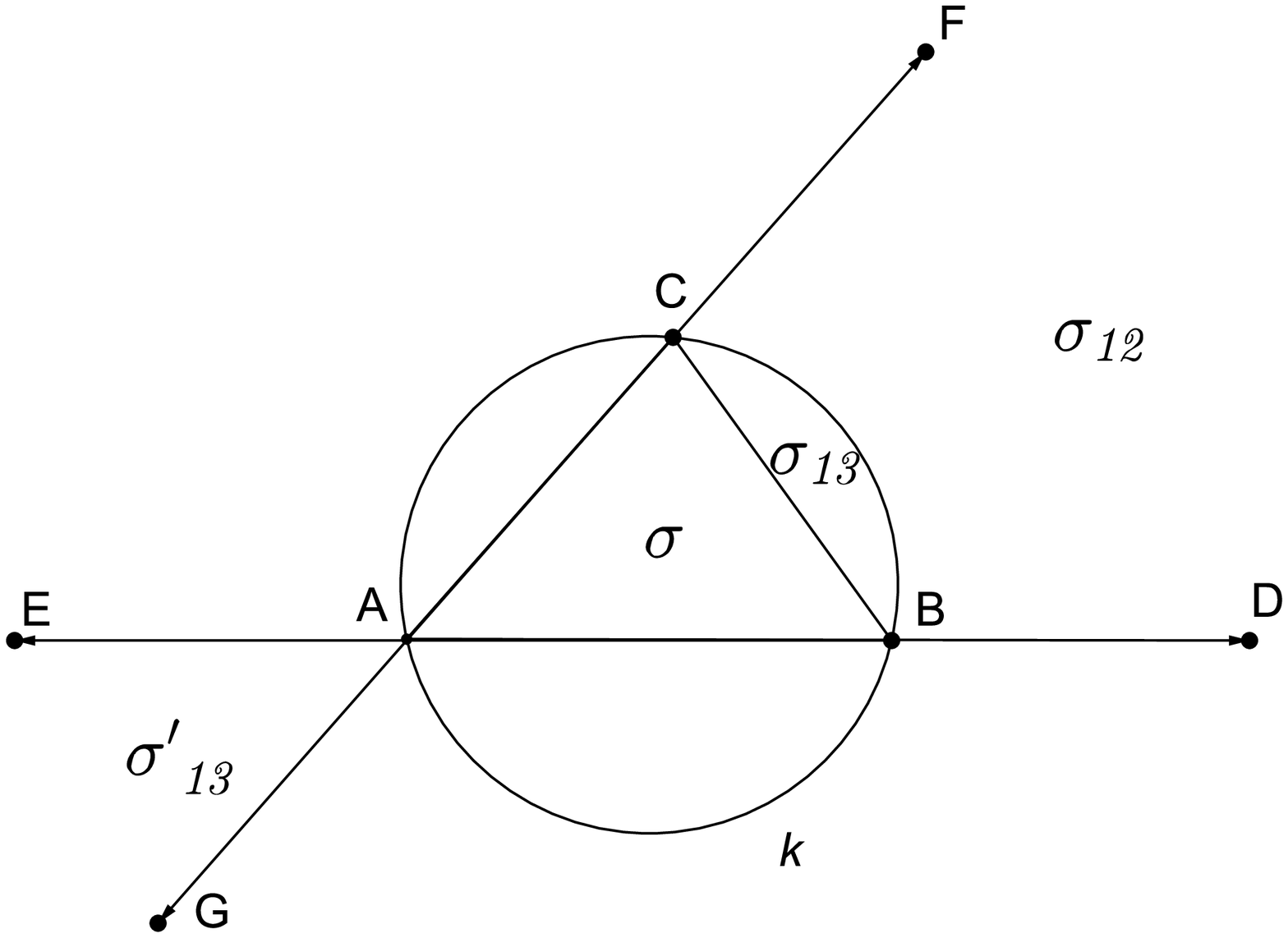}
\vskip 2mm
Fig. 1
\end{center}

1) $\imath(\sigma)=\sigma$.

\noindent
If $M$ is a point on the side $BC$, then $\imath (M)=A$.

2) $\imath(\sigma_{12})=\sigma_{12}$.

\noindent
For any point $M$ on the ray $BD^{\rightarrow}$ we have $\imath (M)=C$; for any point $M \in CF^{\rightarrow}$ $\imath (M)=B$.

3) $\imath(\sigma_{13})=\sigma'_{13}$, \; $\imath(\sigma'_{13})=\sigma_{13}$.

\noindent
For any point $M \in AE^{\rightarrow}$ $\imath(M)=C$; \; for any point $M \in AG^{\rightarrow}$ $\imath(M)=B$.

The transformation $\imath$ can also be defined for points on the circumcircle $k$, different
from the vertices of the triangle. If $M$ is a point on the arc $\widearc {BC}$, then $\imath(M)$ is the point at infinity of the
line, which is symmetric to the line $AM$ with respect to the bisector of $\angle BAC$.

For an arbitrary point $X$ in the plane of
$\triangle ABC$ we denote by $x,\,y,\,z,$ the directed distances of the point $X$
to the lines $BC,\,CA,\,AB$, respectively. Then $(x,\, y,\, z)$ is the triple of trilinear coordinates of $X$
with respect to the basic triangle. The trilinear coordinates satisfy the equality $ax +by +cz = 2S.$
Further we denote by $S_1, \, S_2,\, S_3$ the oriented areas
of the triangles $BCX, CAX, ABX$, respectively. Then
$\ds{\lambda =\frac{S_1}{S}, \, \mu = \frac{S_2}{S}, \, \nu = \frac{S_3}{S}}; \; \lambda + \mu + \nu =1$
are the barycentric coordinates of $X$ with respect to the $\triangle ABC$. The relation between the trilinear coordinates
and the barycentric coordinates of the point $X$ is given by
$$\lambda = \frac{ax}{2S}\,, \quad \mu = \frac{by}{2S}\,, \quad \nu = \frac{cz}{2S}\,.$$

We also consider the function ${\ds J= \frac{a^2}{\lambda}\,+ \frac{b^2}{\mu}\,+
\frac{c^2}{\nu}}=\ds{\left(\frac{a}{x}+\frac{b}{y}+\frac{c}{z}\right)2S,}$ \,
which is defined for all points that do not lie on the sidelines of $\triangle ABC$.

The following characterization of the points in the $\angle BAC$ out of $\triangle ABC$ is useful.
\vskip 2mm
{\bf Lemma 1.} Let $M(x, y, z)$ be with trilinear coordinates satisfying the conditions
$x < 0, \\ y > 0, \, z > 0$, i.e. $M$ is in the $\angle BAC$ out of $\triangle ABC$. Then

1) $M \in \sigma_{12}$ iff $J > 0$;

2) $M$ lies on the arc $\widearc{BC}$ iff $J = 0$;

3) $M \in \sigma_{13}$ iff $J < 0$.
\vskip 2mm
We also need the following statement.

{\bf Lemma 2.} Let $M(x, y, z)$ be with trilinear coordinates satisfying the conditions
$x > 0, \\ y < 0, \, z < 0$, i.e. $M \in \sigma'_{13}$. Then $J<0$.
\vskip 2mm
\textbf{The isogonal conjugation with respect to barycentric coordinates }

Let $(\lambda, \mu, \nu), \, \lambda + \mu +\nu =1$ be the barycentric coordinates of a point $M$,
which does not lie on the lines $AB, BC, CA$ or on the circumcircle $k(ABC)$. If $(\lambda', \mu', \nu')$
are the barycentric coordinates of the point $N$, isogonal conjugate to $M$, then
$$\lambda'=\frac{a^2}{\lambda\, J}\,, \quad \mu'=\frac{b^2}{\mu\, J}\,, \quad \nu'=\frac{c^2}{\nu \,J}\,.\leqno(2.1)$$
If $(\lambda, \mu, \nu)$ are the barycentric coordinates of a point $M$, it is useful to consider the
\emph{homogeneous} barycentric coordinates of $M$:
$$(\rho \lambda,\, \rho \mu,\, \rho \nu), \quad \rho \neq 0.$$

Then the formulas
$$\lambda'=\frac{a^2}{\lambda}\,, \quad \mu'=\frac{b^2}{\mu}\,, \quad \nu'=\frac{c^2}{\nu}\leqno(2.2)$$
represent the isogonal conjugation even on the arcs $\widearc{BC}$, $\widearc{CA}$ or $\widearc{AB}$ of $k$.
If $M(\lambda, \mu, \nu)$ lies on\
the arc $\widearc{BC}$, then the point $\imath(M)= N(\lambda', \mu', \nu')$
satisfies the condition $\lambda'+\mu'+\nu'=0$ and lies on the line at infinity.
\vskip 2mm
\textbf{A general formulation of the problem}
\vskip 2mm
Now, let $(\lambda, \mu, \nu) \neq (0, 0, 0)$ be a triple of fixed real numbers. For any point
$X$ with trilinear coordinates $(x, y, z)$ consider the function
$$F(X)= \lambda x^2+\mu y^2 + \nu z^2, $$
which is a weighted sum of the squares of the directed distances $(x, y, z)$.

Our aim is to investigate the minima and maxima of the above function.

Further we consider three essential cases:

1. $\lambda \mu \nu \neq 0, \quad \lambda + \mu + \nu > 0 $;

2. $\lambda \mu \nu \neq 0, \quad \lambda + \mu + \nu < 0$;

3. $\lambda \mu \nu \neq 0, \quad \lambda + \mu + \nu = 0$.

\section{Weighted sum with $\lambda \mu \nu \neq 0, \quad \lambda + \mu + \nu > 0$}

Obviously both functions
$$\lambda x^2 + \mu y^2 + \nu z^2, \qquad \frac{\lambda x^2 + \mu y^2 + \nu z^2}{\lambda + \mu + \nu}$$
have minima and maxima at the same points. Without loss of generality
we can assume that $\lambda + \mu + \nu = 1$.

Thus the problem in this section is to find the minimum (maximum) of the function
$$F(X)=\lambda x^2 + \mu y^2 + \nu z^2, \quad \lambda + \mu + \nu = 1. \leqno(3.1)$$

First we consider the case
\vskip 2mm
\textbf{1.1.} $\lambda > 0, \; \mu > 0, \; \nu > 0$.
\vskip 2mm
\textbf{Problem 1.} {\it
$(i)$ Find the point $N$ that minimizes the function $F(X)$.

$(ii)$ If $M$ is the point with barycentric coordinates $(\lambda, \mu, \nu)$, prove that $M$ and $N$ are
isogonal conjugate.}

{\it Solution.} To solve (i), consider the system
$$\begin{array}{l}
{\ds F(X)=\lambda\,x^2+\mu\,y^2+\nu\,z^2,}\\
[2mm]
ax+by+cz=2S;
\end{array}
\qquad
x, y , z \in \R \leqno(3.2)
$$
and interpret $(x,y,z)$ as Cartesian coordinates in the three dimensional Euclidean space.
The level surfaces of the function $F(X)$ are the ellipsoids
$$\varepsilon(k) : \quad
\lambda\,x^2+\mu\,y^2+\nu\,z^2=k, \quad k={\rm const} \in (0, \infty).$$

Geometrically, to find the point that minimizes the function $F(X)$ in (2.2), means to find $k$ so that the ellipsoid $\varepsilon(k)$
is tangent to the plane $\pi: ax+by+cz=2S$ and then to determine the touch-point $N$ of $\varepsilon(k)$ to $\pi$.

The tangent plane to $\varepsilon(k)$ at a point $(x_0, y_0, z_0)$ is given by the
equality
$$\tau: \lambda\,x_0 \, x+\mu\,y_0 \, y+\nu\,z_0 \, z=k.$$

Then the condition $\tau \equiv \pi$ implies that
$$\begin{array}{l}
\ds{\frac{\lambda\,x_0}{a}=\frac{\mu\,y_0}{b}=\frac{\nu z_0}{c}=t,}\\
[4mm]
t(a x_0+b y_0+c z_0)=k,\\
[4mm]
ax_0+by_0+cz_0=2S.
\end{array} \leqno(3.3)
$$

Solving (3.3), we find:
$$t= \frac{2S}{J}\,, \quad k=\frac{4S^2}{J}\,,\; \left(J=\frac{a^2}{\lambda}+
\frac{b^2}{\mu}+\frac{c^2}{\nu}\right)$$
and
$$x_0=\frac{2S}{J}\,\frac{a}{\lambda}\,, \quad
y_0=\frac{2S}{J}\,\frac{b}{\mu}\,, \quad z_0=\frac{2S}{J}\,\frac{c}{\nu}\,.\leqno(3.4)$$

Now, taking into account (3.4), we conclude that the point $N$ minimizing the function $F(X)$
has barycentric coordinates $(\lambda', \mu', \nu')$ with respect to $\triangle ABC$
given by
$$\lambda'=\frac{ax_0}{2S} = \frac{a^2}{\lambda\,J}\,, \quad \mu'=\frac{b\,y_0}{2S} = \frac{b^2}{\mu\,J}\,,
\quad \nu'=\frac{cz_0}{2S} = \frac{c^2}{\nu \, J}\, \leqno(3.5)$$
and $\ds{F_{min} = k = \frac {4S^2}{J}}$, which solves (i).

To prove (ii), let us denote by $M$ the point with barycentric coordinates $(\lambda, \mu, \nu)$.
Comparing with (1.1) we conclude that formulas (3.5) are a representation of the isogonal conjugation in
barycentric coordinates. Hence, the point $N(\lambda', \mu', \nu')$ is the isogonal conjugate
one to the point $M(\lambda, \mu, \nu)$.
\qed

In this case $M$ and $N$ are in $\sigma$.

Next we consider the case
\vskip 2mm
\textbf{1.2.} $\lambda < 0, \; \mu > 0, \; \nu > 0$.
\vskip 2mm
In this case the problem states as follows:

\textbf{Problem 2.} {Prove that $F(X)$ has a minimum if and only if
$$\ds{J=\frac{a^2}{\lambda}+\frac{b^2}{\mu}+\frac{c^2}{\nu} < 0}.$$

$(i)$ Find the point $N$ that minimizes the function $F(X)$.

$(ii)$ If $M$ is the point with homogeneous barycentric coordinates $(\lambda, \mu, \nu)$, prove that $M$ and $N$ are
isogonal conjugate.}

{\it Solution.} Consider the system (3.2). In this case any level surface of the function $F(X)$
$$\varepsilon(k) : \quad
\lambda\,x^2+\mu\,y^2+\nu\,z^2=k, \quad k={\rm const} \in \R$$
is one of the following: one sheet hyperboloid if $k > 0$; cone if $k = 0$; two sheet hyperboloid if $k < 0$.

Geometrically, to find the point that minimizes the function $F(X)$ in (3.2), means to find $k<0$ so that the two sheet
hyperboloid $\varepsilon(k)$ is tangent to the plane $\pi: ax+by+cz=2S$ and then to determine the touch-point
$N$ of $\varepsilon(k)$ to $\pi$.

The plane $\pi$ can be tangent to $\varepsilon(k)$ only if $ k< 0 $.

Similarly to the solution of Problem 1 we obtain the system (3.4), which implies that
$$t J= 2S, \quad 2S t = k. $$
Therefore $\varepsilon(k)$ is tangent to the plane $\pi$ only in the case $J < 0$.

Further, we find the coordinates of the touch-point
$$x_0=\frac{2S}{J}\,\frac{a}{\lambda}\,, \quad
y_0=\frac{2S}{J}\,\frac{b}{\mu}\,, \quad z_0=\frac{2S}{J}\,\frac{c}{\nu}\,.\leqno(3.6)$$

Thus, the point $N$ minimizing the function $F(X)$ lies in the domain  $\sigma'_{13}$.

Let $M$ be the point with barycentric coordinates $(\lambda, \, \mu, \, \nu)$. Then the formulas (3.6) show that
the point $N(\lambda', \, \mu', \, \nu')$ has barycentric coordinates
$$\lambda'=\frac{a^2}{\lambda\,J}\,, \quad \mu'=\frac{b^2}{\mu\,J}\,, \quad \nu'=\frac{c^2}{\nu\,J}$$
and it is isogonal conjugate to the point $M$.

Hence, the triple $(\lambda, \mu, \nu)$ determines a point $M$ in the domain $\sigma_{13}$ and
$$F_{min}= F(N) = \frac{4S^2}{J} < 0,$$
where $N \in \sigma'_{13}$ is the isogonal conjugate to the point $M$.
\qed
\vskip 2mm
\textbf{1.3.} $\lambda >0,\; \mu < 0,\; \nu < 0$.

In this case the problem states as follows:

\textbf{Problem 3.} {Prove that the function $F(X)$ has neither minimum, nor maximum.}

{\it Solution.} Let $M$ be the point with barycentric coordinates $(\lambda, \, \mu, \, \nu)$.
Then $M \in \sigma'_{13}$ and $J < 0$. Similarly to the case 1.2 we obtain that
in our case any level surface of the function $F(X)$
$$\varepsilon(k) : \quad
\lambda\,x^2+\mu\,y^2+\nu\,z^2=k, \quad k={\rm const} \in \R$$
is one of the following: one sheet hyperboloid if $k < 0$; cone if $k = 0$; two sheet hyperboloid if $k > 0$.
Since $J < 0$, then $\ds{k=\frac{4S^2}{J} < 0}$ and $\pi$ can not be a tangent plane to any two sheet hyperboloid.
Hence the function $F(X)$ has neither minimum nor maximum.
\qed
\vskip 2mm
\section{Weighted sum with $\lambda \mu \nu \neq 0, \; \lambda + \mu + \nu < 0$}
If $\lambda + \mu + \nu < 0$, then we put $\bar \lambda = - \lambda,\; \bar \mu = - \mu, \; \bar \nu = - \nu$
and consider the function
$$\bar F(X) = \bar \lambda x^2 + \bar \mu y^2 + \bar \nu z^2=- F(X),  \quad \bar\lambda + \bar\mu + \bar\nu > 0.$$
We suppose again that $\bar\lambda + \bar\mu + \bar\nu = 1$ and consider the point $M$ with barycentric coordinates
$(\bar\lambda, \bar\mu, \bar\nu )$.

Comparing with Section 2 we have the following.
\vskip 2mm
\textbf{2.1. $\lambda<0, \; \mu <0, \; \nu < 0$}

Under these conditions $\ds{J=\frac{a^2}{\lambda}+\frac{b^2}{\mu}+\frac{c^2}{\nu} < 0}.$

Then
$$M(\bar \lambda, \bar \mu, \bar \nu) \in \sigma, \; N=\imath(M) \in \sigma, \; \bar J>0, \quad
F_{max}=F(N)= \frac{4S^2}{J} < 0.$$

\textbf{2.2. $\lambda>0, \, \mu <0, \; \nu < 0$}

\textbf{2.2.1. $J<0$}

The function $\bar F =-F$ has no minimum or maximum.

\textbf{2.2.2 $J>0$}
$$M(\bar \lambda, \bar \mu, \bar \nu) \in \sigma_{13}, \; N=\imath(M) \in \sigma'_{13}, \; \bar J<0, \quad
F_{max}=F(N)= \frac{4S^2}{J} > 0.$$

\textbf{2.3. $\lambda<0, \; \mu >0, \; \nu > 0$}

The function $\bar F =-F$ has no minimum or maximum.

\section{Weighted sum with $\lambda \mu \nu \neq 0, \; \lambda + \mu + \nu = 0$}

Let us consider the system of mass points $\{A(\lambda), B(\mu), C(\nu)\}$ and denote by $P(\mu+\nu)$ the center of
mass of the system $\{B(\mu), C(\nu)\}$ (Fig. 2). If $O$ is an arbitrary point, we have
$$\lambda \overrightarrow{OA} + \mu \overrightarrow{OB} + \nu \overrightarrow{OC}=-\lambda \overrightarrow{AP}
=\mu \overrightarrow{AB} + \nu \overrightarrow{AC} = \overrightarrow{v} = \overrightarrow{\rm const}.$$
Calculating
$$\overrightarrow{v}^2=-\lambda \mu \nu \left(\frac{a^2}{\lambda}+\frac{b^2}{\mu}+\frac{c^2}{\nu}\right)=
- \lambda \mu \nu J,$$
we obtain the geometric meaning of $J$ in the case $\lambda + \mu + \nu = 0$.

The condition $\overrightarrow{v}^2> 0$ implies that
$$\lambda \mu \nu J < 0. \leqno(5.1)$$

\begin{center}
\includegraphics[width=10cm]{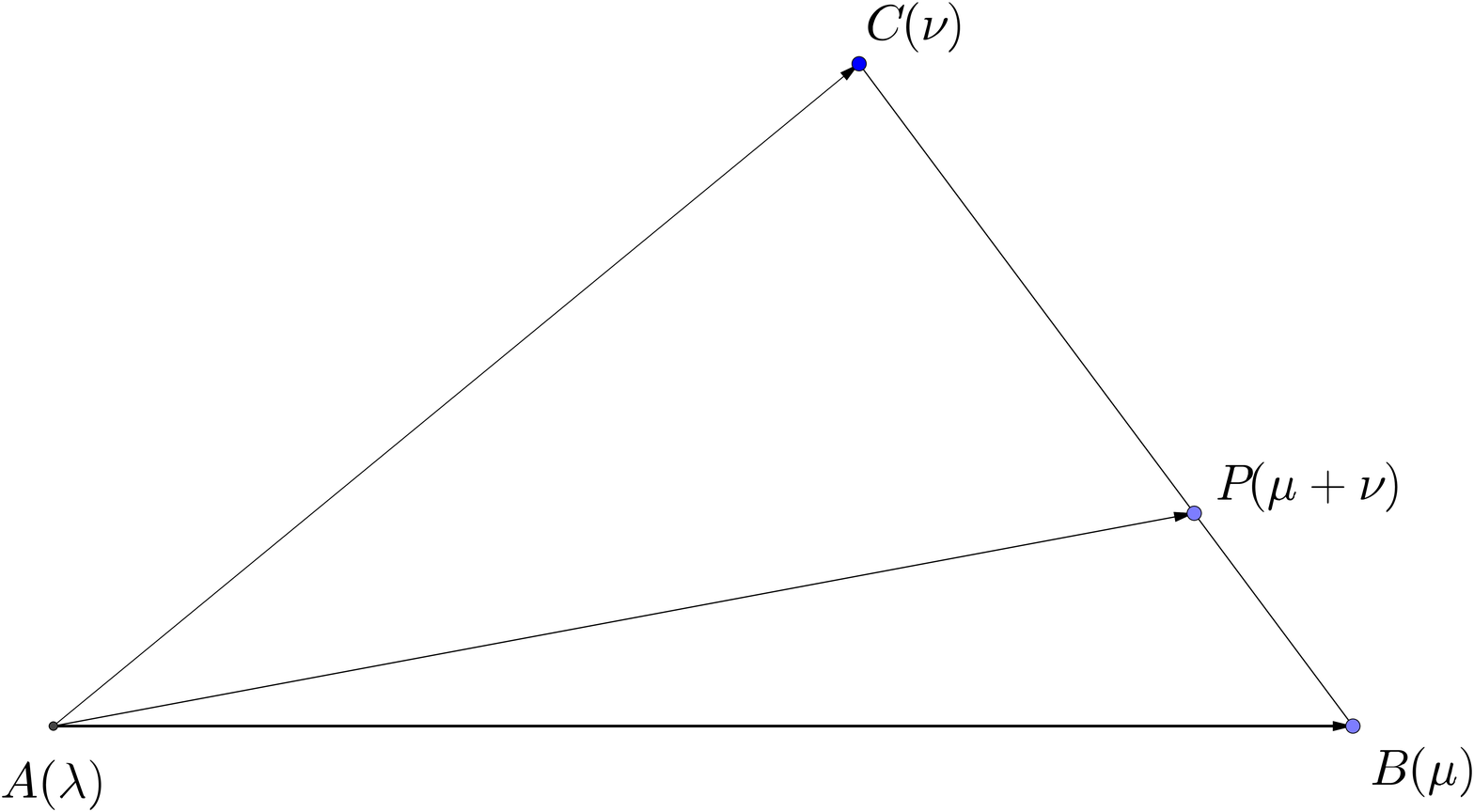}
\vskip 2mm
Fig. 2
\end{center}

In this section we consider two cases: $(\lambda <0, \mu > 0, \nu > 0)$ and $(\lambda > 0, \mu < 0, \nu < 0)$.

\textbf{3.1. $\lambda \mu \nu < 0.$}

In this case $M(\lambda, \mu, \nu)$ can be interpreted as a point at infinity, i.e. the point at infinity of
the line $AP$ (Fig. 2). The inequality (5.1) implies that
$J>0$. Using similar arguments as in Section 2 we conclude that the function $F(X)$ has neither minimum
nor maximum.

\textbf{3.2. $\lambda \mu \nu > 0.$}

Under these conditions the inequality (5.1) implies that $J<0$.
We consider the function $\bar F=-F$. Since $\bar \lambda \bar \mu \bar \nu <0$, comparing with the previous case
we conclude that $\bar F$ i.e. $F$ has no minimum or maximum.

\section{Weighted sum with one or two zero weights}
\vskip 2mm
\textbf{4.1.} $\lambda = 0, \, \mu > 0, \, \nu > 0.$
\vskip 2mm
In this case the level surfaces of the function $F(X)$ are the elliptic cylinders
$$\mu\,y^2+\nu\,z^2=k, \quad k={\rm const} \in (0, \infty),$$
and the axis $Ox$, when $k=0$. The level surfaces intersect the plane $ax +by + cz = 2S$ into ellipses
by $k > 0$ and in the point $(2S, 0, 0)$ by $k=0$.

The minimum of the function $F(X)$ is $F_{min}=0$ and it occurs when $y=z=0$.

Thus $M(0, \mu, \nu)$ is a point on the sideline $BC$ and the minimum $F_{min}=0$
occurs when $N\equiv A$, i.e. $M$ and $N$ are again isogonal conjugate.
\vskip 2mm
\textbf{4.2. $\lambda=0, \; \mu <0, \; \nu < 0$.}
\vskip 2mm
$M(0, \bar \mu, \bar \nu) \in BC, N=\imath(M)=A,$ $F_{max}=F(N)=0$.
\vskip 2mm
\textbf{4.3. $\lambda=0, \; \mu \nu <0 $.}
\vskip 2mm
In this case the level surfaces of the function $F(X)$
$$\mu y^2 + \nu z^2 = k, \quad k\in (-\infty,+\infty)$$
are hyperbolic cylinders if $k \neq 0$ and two planes if $k=0$.

Hence the function $F(X)$ has no minimum or maximum.
\vskip 2mm
\textbf{5.1.} $\lambda = 1, \, \mu = \nu = 0.$
\vskip 2mm
It is clear that $F_{min}=0$ and it occurs when $x=0$.

Thus $M(1, 0, 0) \equiv A$ and $N $ is any point on $BC$, i.e. $M$ and $N$ are again isogonal conjugate.
\vskip 2mm
\textbf{5.2. $\lambda=-1, \; \mu =0, \; \nu = 0$.}
\vskip 2mm
$M(1,0,0)\equiv A, N \in BC$, $F_{max}=F(N)=0$.
\vskip 2mm
\textbf{Summarizing we get the following:}
\vskip 2mm
1. $\lambda + \mu + \nu > 0$.

If $M(\lambda, \mu, \nu)$ is an inner point for the circumcircle $k$, or coincides
with a vertex of $\triangle ABC$, then $F$ has a minimum and $F_{min}=F(N), \, N=\imath(M)$.
\vskip 2mm
2. $\lambda + \mu + \nu < 0$.

If $M(-\lambda, -\mu, -\nu)$ is an inner point for the circumcircle $k$, or coincides
with a vertex of $\triangle ABC$, then $F$ has a maximum and $F_{max}=F(N), \, N=\imath(M)$.
\vskip 2mm
3. $\lambda + \mu + \nu = 0$.

In this case the function $F$ has neither a minimum nor a maximum.

\section{Examples}

We choose as a typical example the following pair of conjugate points: the circumcenter $O$ and the orthocenter $H$.

1. Let $M\equiv O$. Then $\lambda =\sin 2\alpha, \, \mu = \sin 2\beta, \, \nu = \sin 2 \gamma$ and the
point $O$ generates the function
$$F(X)=\sin 2\alpha \, x^2+\sin 2\beta \, y^2+\sin 2 \gamma \, z^2.$$

1.1. $\triangle ABC$ is acute-angled and $O \in \sigma$. Simple calculations show that
$$J=\frac{4S}{\cos \alpha \, \cos \beta \, \cos \gamma}, \quad F_{min}=F(H)= 4S \cos \alpha \, \cos \beta \, \cos \gamma.$$

1.2. $\alpha = 90^0$ and $O$ is the midpoint of $BC$. Then
$$F_{min}=F(A)=0.$$

1.3. $\alpha > 90^0$ and $O \in \sigma_{13}$. Then
$$J=\frac{4S}{\cos \alpha \, \cos \beta \, \cos \gamma}\,<0, \quad
F_{min}=F(H)= 4S \cos \alpha \, \cos \beta \, \cos \gamma \, <0.$$
\vskip 2mm

{\bf Remark.} \emph{The function F(X) generates geometric inequalities.}

Let $M(\lambda, \, \mu, \, \nu)\in \sigma$ generate the function (3.1). If we choose a concrete triangle
center $X$ in $\sigma$ with trilinear coordinates $(x, y, z)$ and replace in (3.1), then we obtain
the geometric inequality
$$F(X)\geq F_{min}=\frac{4S^2}{J}\,.$$
The equality occurs if and only if $X\equiv \imath(M)$.

In a similar way the case $M \in \sigma_{13}$ generates even more interesting geometric inequalities.

\end{document}